\documentclass{amsart}
\usepackage{amsmath}
\usepackage{graphicx}
\usepackage{latexsym}
\usepackage{color}
\usepackage{amscd}
\usepackage[all]{xy}
\usepackage{enumerate}
\usepackage{hyperref}
\theoremstyle{plain}
\newtheorem{theorem}{Theorem}[section]

\newtheorem{lemma}[theorem]{Lemma}

\theoremstyle{definition}

\newtheorem{remark}[theorem]{Remark}

\newtheorem*{mthm}{Theorem~A}
\newtheorem*{mcor}{Theorem~B}



\newcommand{\CPb}{\overline{\mathbb{CP}}{}^{2}}
\newcommand{\CP}{{\mathbb{CP}}{}^{2}}

\newcommand{\N}{\mathbb{N}}
\newcommand{\R}{\mathbb{R}}

\newcommand{\Z}{\mathbb{Z}}

\newcommand{\K}{{\rm K3}}

\def \x {\times}

\begin{document}

\title[Inequivalent Lefschetz fibrations and surgery equivalence]
{Inequivalent Lefschetz fibrations and \\ surgery equivalence of symplectic $4$-manifolds}

\author[R. \.{I}. Baykur]{R. \.{I}nan\c{c} Baykur}
\address{Department of Mathematics and Statistics, University of Massachusetts, Amherst, MA 01003-9305, USA}
\email{baykur@math.umass.edu}

\begin{abstract}
We prove that any symplectic $4$-manifold which is not a rational or ruled surface, after sufficiently many blow-ups, admits an arbitrary number of nonisomorphic Lefschetz fibrations of the same genus which cannot be obtained from one another via Luttinger surgeries. This generalizes results of Park and Yun who constructed pairs of nonisomorphic Lefschetz fibrations on knot surgered elliptic surfaces. In turn, we prove that there are monodromy factorizations of Lefschetz pencils which have the same characteristic numbers but cannot be obtained from each other via partial conjugations by Dehn twists, answering a problem posed by Auroux. 
\end{abstract}

\maketitle

\section{Introduction} 

Since Donaldson's proof of the existence of Lefschetz pencils on symplectic \linebreak $4$-manifolds \cite{Donaldson}, an immense literature has been dedicated to the study of Lefschetz fibrations and pencils. However, a rather fundamental question has been left mostly unanswered: 

\indent \textit{How many nonisomorphic Lefschetz pencils\,/\,fibrations of the same genus does} \linebreak 
\indent \textit{a given symplectic $4$-manifold admit?} \cite{ParkYun} 

\noindent Here, two Lefschetz pencils\,/\,fibrations are called  \emph{isomorphic} if there are orientation-preserving self-diffeomorphisms of the $4$-manifold and the base surface which make the two fibrations commute ---where it is clear that the fiber genera, as well as the number of base points in the case of pencils, should match to begin with. 

Complementary to the diversity question above, one can inquire about how to relate two different Lefschetz pencils\,/\,fibrations. Let $(X,f)$ be a symplectic Lefschetz pencil and $L \subset X \setminus \text{Crit(f)}$ be an embedded Lagrangian torus that fibers over a loop $\gamma$ on the base, i.e. $f|_L$ is a circle bundle over $\gamma$ obtained by a parallel transport of a loop $\alpha$ in the fiber. A \textit{Luttinger surgery} along $L$ in the direction of $\alpha$, which we will call a \textit{fibered Luttinger surgery} in short, yields a new symplectic Lefschetz pencil, which has the same Euler characteristic, signature and symplectic Kodaira dimension \cite{Li3} as $X$, and if $(X,f)$ is supported by an \textit{integral} symplectic form $\omega$, the characteristic numbers $[\omega]^2$ and $c_1\cdot[\omega]$ do not change \cite{ADK}. Here is the second question we are interested in:

\indent \textit{Are any pair of integral symplectic Lefschetz pencils with the same characteristic} \linebreak
\indent \textit{numbers $(c_1^2, c_2, [\omega]^2$, $c_1 \cdot [\omega])$ related via (fibered) Luttinger surgeries?}  \cite{Auroux2}

\smallskip
\noindent This is a Lefschetz pencil version of the --still open-- question on the equivalence of integral symplectic $4$-manifolds with the same characteristic numbers listed above via Luttinger surgeries \cite{Auroux1}, which in turn is a symplectic version of Stern's question on the equivalence of homeomorphic smooth $4$-manifolds via smooth surgeries along tori \cite{Stern}, settled positively in \cite{BS}.

Our goal is to address both questions above by proving the following:  

\begin{mthm} \label{mainthm}
Given any closed symplectic $4$-manifold $X$ which is not a rational or ruled surface, and any positive integer $N$, there are $N$ nonisomorphic Lefschetz pencils of the same genus on a blow-up of $X$, which are not equivalent via fibered Luttinger surgeries. These pencils can be chosen so that they only have nonseparating vanishing cycles.
\end{mthm}

In fact, we will show in the proof of the theorem that there is no sequence of --not necessarily fibered-- Luttinger surgeries starting from the blow-up of $X$ and leading back to it so that the fiber of one pencil transforms into the fiber of another. Thus the rigidity here is not due to the choice of Lagrangian tori and surgery coefficients, but due to relating two pencils instead of underlying symplectic manifold(s). As usual, blowing-up all the base points, one obtains inequivalent Lefschetz fibrations with the same properties listed in the theorem.

For a quick insight into the correlation between the above questions of seemingly different nature, recall that an isomorphism between two Lefschetz fibrations can be translated to a combinatorial equivalence between the associated positive Dehn twist factorizations in the mapping class groups of surfaces via \textit{Hurwitz moves} and \textit{global conjugations} \cite{Kas, Matsumoto}. A weaker equivalence between two Lefschetz fibrations on two possibly different symplectic $4$-manifolds can be then defined by allowing \textit{partial conjugations}, i.e. conjugations of positive subwords of the monodromy factorizations in the mapping class group. A particularly important case of partial conjugation is when the conjugating mapping class is a Dehn twist along a curve $\alpha$ preserved by the positive subword, which corresponds to a fibered Luttinger surgery as above \cite{Auroux2}. 

In \cite{ParkYun}, Park and Yun appealed to this combinatorial approach to show that there are pairs of inequivalent Lefschetz \textit{fibrations} on certain knot surgered elliptic surfaces produced by Fintushel and Stern (and it is mentioned in \cite{ParkYun} that Smith had another isolated example of pairwise distinct Lefschetz fibrations on $T^2 \x \Sigma_2 \# 9 \CPb$ in his thesis). In order to obtain their result, the authors proved that for certain pairs of fibered knots, one gets two distinct subgroups of the mapping class group generated by the \textit{collection} of positive Dehn twists in respective monodromy factorizations. Curiously, all these examples were on symplectic $4$-manifolds of Kodaira dimension $1$, did not contain $(-1)$-sphere sections (i.e. they were not blow-ups of pencils), and were indeed equivalent via partial conjugations. Moreover, Park and Yun's clever use of monodromy groups in \cite{ParkYun} was not powerful enough to distinguish more than two fibrations. Our theorem generalizes their result in all these aspects, generating examples on blow-ups of almost all symplectic $4$-manifolds. 

On the other hand, Auroux posed the following \cite[Question 5]{Auroux2}: 

\indent \textit{Are any two factorizations of a boundary parallel twist into positive Dehn }\linebreak
\indent \textit{twists along nonseparating curves always equivalent via Hurwitz moves and } \linebreak
\indent \textit{partial conjugations, provided the two factorizations give the same Euler } \linebreak
\indent \textit{characteristic and signature?} 

\smallskip
\noindent As implicitly suggested by the question, one can read off the Euler characteristic and signature of the total space of such a fibration from the factorization (where the latter, given by the Meyer signature cocycle, is a much harder calculation, but possible in principle). Furthermore, since the problem is formulated in analogy with the Luttinger surgery question above, the particular interest here is in the case when the partial conjugation along a subword is performed by a Dehn twist along a curve that is fixed with the same orientation by this subword. We will call this an \textit{untwisted partial conjugation}. Our second theorem, which is a weaker reformulation of Theorem~A in this setting, answers the above question in the negative:

\begin{mcor} \label{maincor}
The positive boundary multitwist in the mapping class group $\Gamma_g^m$ of a closed orientable genus $g$ surface with $m$ boundary components admits $N$ factorizations into positive Dehn twists along nonseparating curves, which have the same Euler characteristic and signature, but are not equivalent via Hurwitz moves and untwisted partial conjugations by Dehn twists, where $g$ is taken sufficiently large for given $N, m \geq 1$. 
\end{mcor}

Our proof of Theorem~A does not deal with explicit monodromy factorizations, and instead uses a variant of the degree doubling construction from \cite{Smith, AK} for \textit{topological} Lefschetz pencils. In Section~\ref{Proof} we are going to introduce an algebraic invariant defined on equivalence classes of Lefschetz pencils up to isomorphisms and (fibered) Luttinger surgeries, using arguments that rely on Seiberg-Witten theory \cite{Ta,Li1} and holomorphic curves \cite{Welschinger}. This Lefschetz pencil invariant, which we call \emph{exceptional data}, encodes the number of certain multisections representing exceptional classes in the symplectic $4$-manifold. Its importance to us rests on the easiness in keeping track of how exceptional data changes under \emph{partial doubling sequences} we will introduce, which simply consist of blow-ups/blow-downs and degree doublings for pencils, and in turn, this will allow us to reduce the problem to a merely combinatorial one. All these are contained in Section~\ref{Proof}, where we prove Theorem~A. In Section~\ref{Partial} we prove Theorem~B and discuss related problems on cobordisms between symplectic $4$-manifolds.

\vspace{0.3in}
\noindent \textit{Acknowledgements.} The author would like to thank Denis Auroux and Jeremy Van Horn-Morris for their helpful comments on a draft of this paper, and the anonymous referee for kindly pointing out an error in the proof of Theorem~B in the earlier version. The author was partially supported by the NSF Grant DMS-$1510395$, the Simons Foundation Grant $317732$ and the ERC Grant LDTBud.

\section{Background} \label{Background}

Here we review the definitions and basic properties of Lefschetz pencils and fibrations, mapping class groups, and Luttinger surgeries. The reader can turn to \cite{GS} and \cite{ADK} for more details.

Let $X$ be a closed, oriented $4$-manifold, and $B=\{b_j \}$, $C=\{p_i \}$ be finite \textit{non-empty} sets of points in $X$. A \emph{Lefschetz pencil} $(X,f)$ is a surjection $f$ from $X \setminus B$ onto $S^2$ that is a submersion on the complement of $C$, such that around each \emph{base point} $b_j$ and \emph{critical point} $p_i$ there are local complex coordinates (compatible with the orientations on $X$ and $S^2$) with respect to which the map takes the form \linebreak $(z_1,z_2) \mapsto z_1/z_2$ and $(z_1, z_2) \mapsto z_1 z_2$, respectively.  A \emph{Lefschetz fibration} is defined similarly when $B = \emptyset$. Blowing-up all the base points $b_j$ in a pencil $(X,f)$, one obtains a Lefschetz fibration $(X',f')$ with disjoint $(-1)$-sphere sections $S_j$ corresponding to $b_j$. We say that we have a \emph{genus $g$ pencil/fibration}, for $g$ the genus of the regular fiber $F$ of the pencil/fibration (which is compactified by including the base points). The fiber containing the critical point $p_i$ has a nodal singularity at $p_i$, which locally arises from shrinking a simple loop $a_i$ on $F$, called the \emph{vanishing cycle}. A singular fiber of a Lefschetz pencil/fibration $(X,f)$ is called \emph{reducible} (resp. \emph{irreducible}) if $a_i$ is separating (resp. nonseparating). 

By the seminal work of Donaldson every symplectic $4$-manifold $(X, \omega)$ admits a \emph{symplectic} Lefschetz pencil whose fibers are symplectic with respect to $\omega$ \cite{Donaldson}. Conversely, Gompf showed that total space of a Lefschetz fibration (recall $C \neq \emptyset$), and in particular blow-up of any pencil, always admits a symplectic form $\omega$ with respect to which all regular fibers and any chosen collection of disjoint sections are symplectic, and in fact, any such two symplectic forms are deformation equivalent \cite{GS}. We will often use the notation $(X, \omega, f)$ to indicate that $f$ is a symplectic Lefschetz pencil/fibration with respect to $\omega$, where any explicitly discussed sections of $f$ will always be assumed to be symplectic as well.  

Lefschetz pencils and fibrations can be described combinatorially in terms of products of Dehn twists in the mapping class group as follows: Let $\Sigma_{g}^m$ denote a compact oriented surface of genus $g$ with $m$ boundary components, with the convention that $\Sigma_g=\Sigma_{g}^0$. The \emph{mapping class group\,}, $\Gamma_{g}^m$, of  $\Sigma_{g}^m$ is the group of orientation-preserving self-diffeomorphisms of $\Sigma_{g}^m$ fixing the points on the boundary up to isotopies fixing the points on the boundary as well. Let $t_a \in \Gamma_g^m$ denote the positive (right-handed) Dehn twist along the simple loop $a$ on $\Sigma_g^m$. Now let us also assume that all critical points $p_i$ of the Lefschetz pencil/fibration $(X,f)$ lie in distinct fibers, which can always be achieved after a small perturbation. It turns out that the local monodromy around the singular fiber with vanishing cycle $a_i$ is $t_{a_i}$, and thus, the global monodromy of the fibration around all the singular fibers (i.e. in the complement of a regular fiber) is a product $t_{a_1} \cdot \ldots \cdot t_{a_r}$ in $\Gamma_g^m$, called the \emph{monodromy factorization}, where $r=|C|$. The fact that the map extends to the neighborhood of the excluded regular fiber dictates that the relation
\[ t_{a_1} \cdot \ldots \cdot t_{a_r} = t_{\delta_1} \cdot \ldots \cdot t_{\delta_m} \]
holds in $\Gamma_g^m$, where $\delta_j$ are boundary components of $\Sigma_g^m$. Conversely, provided $g \geq 2$, from any relation between Dehn twists in $\Gamma_g^m$ as above, one can construct a genus $g$ Lefschetz pencil (resp. fibration) with $m$ base points (resp. $m$ disjoint $(-1)$-sphere sections), and $r$ critical points with vanishing cycles $a_i$. Noting all the choices involved in this correspondence, and assuming $g \geq 2$, we get a one-to-one equivalence between Lefschetz pencils up to isomorphisms (i.e. orientation-preserving self-diffeomorphisms of the $4$-manifold and the base surface which make the fibrations commute) and monodromy factorizations up to \textit{Hurwitz moves} and \textit{global conjugations} (i.e. trading subwords $t_{a_i}t_{a_{i+1}}$ with $t_{a_{i+1}}t_{a_{i+1}}^{-1}t_{a_i}t_{a_{i+1}}= t_{a_{i+1}}t_{t_{a_{i+1}}(a_i)}$, and \emph{every} $t_{a_i}$ with $t_{\phi{a_i}}$, $\phi \in \Gamma_g^n$, respectively) \cite{Kas, Matsumoto}. 

Lastly, we review the surgery along Lagrangian tori we are interested in here:\linebreak Let $L$ be an embedded torus in $X$ with trivial normal bundle, identified as\linebreak $\nu L \cong T^2 \x D^2$ under a chosen framing. Given a simple loop $l$ on $L$, let $S^1_{l}$ be a loop on the boundary $\partial(\nu L)\cong T^3$ that is parallel to $l$ under the chosen framing.  Let $\mu_{L}$ denote a meridian circle to $L$ in $\partial(\nu L)$. The $p/q$\/ surgery on $L$ with respect to $l$, describes the smooth $4$-manifold 
\begin{equation*}
X_{L,l,p/q}= (X\setminus\nu L)
\cup_{\varphi} (T^2\times D^2),
\end{equation*}
where the gluing diffeomorphism $\varphi:T^2\times\partial D^2\rightarrow
\partial(X\setminus\nu L)$ satisfies
\begin{equation*}
\varphi_{\ast}([\partial D^2])= p[\mu_{L}] + q[S^1_{l}]
\in H_1(\partial(X\setminus\nu L);\Z).
\end{equation*}
\textit{Luttinger surgery} \cite{Luttinger} is a symplectic analogue of this construction: Let $L$ be an embedded Lagrangian torus in the symplectic $4$-manifold $(X, \omega)$, it then has a canonical framing for $\nu L \cong T^2\times D^2$, called the {\it Lagrangian framing}, such that $T^2\times\{ x \}$ corresponds to a Lagrangian submanifold of $X$\/ for every $x\in D^2$. Using this framing in the above construction, the \emph{Luttinger surgery} on $(X, \omega)$ is a $1/q$ surgery along $L$ with respect to $l$, producing a new $4$-manifold $X' = X_{L,l,1/q}$ with a symplectic form $\omega'$ that restricts to $\omega$ in the complement of the surgery region \cite{ADK}. 

\textit{Fibered Luttinger surgery} along $L$ is then described via choices that count in the fibration structure as described in the Introduction, and amounts to a new monodromy factorization obtained by an untwisted partial conjugation of a subword by a Dehn twist.

\vspace{0.1in}
\section{Proof of Theorem~A} \label{Proof}

The first ingredient we need is the ``degree doubling'' procedure which produces a new genus $g'$ symplectic Lefschetz pencil $(X,\omega, f')$ with $m'$ base points from a given genus $g$ symplectic Lefschetz pencil $(X,\omega, f)$ with $m$ base points, where $g'=2g+m-1$ and $m'=4m$. This construction is described for holomorphic pencils as well as for Donaldson's pencils in Smith's work \cite{Smith}, and for pencils obtained via branched coverings of $\CP$ by Auroux and Katzarkov in \cite{AK} with an explicit calculation of the monodromies. As Smith shows, given a Donaldson pencil with fiber class Poincar\'{e} dual to $d \, [\omega]$, one can pass to a pencil with fiber class dual to $2d \, [\omega]$ using this construction, where the latter pencil has only \textit{nonseparating} vanishing cycles \linebreak \cite[Theorem~3.10]{Smith}. We will use a slight variation of this doubling procedure repeatedly: Let $(X, \omega, f)$ be a Donaldson type symplectic Lefschetz pencil with $m  \geq 1$ base points. Then its \emph{partial double along $m \geq k \geq 1$ points} is the Lefschetz pencil one gets by \textit{first} symplectically blowing-up $(X, \omega, f)$ at $m-k$ points and then doubling the resulting pencil on $(X', \omega', f')$, where $X' = X \# (m-k) \CPb$.

We will need to grant that we can take the double of a given (topological) pencil $(X,f)$ in general. As evident in \cite{AK}, the double $(X,f')$ of the pencil $(X,f)$ is obtained by gluing two pieces; the ``convex'' piece, which is the complement of the regular fiber, is contained in the new pencil as a subpiece, whereas the ``concave'' piece is obtained by a standard doubling of the regular neighborhood of the fiber which is a standard symplectic disk bundle of degree $m$ over a genus-$g$ surface. Thus the \textit{universality property} discussed in \cite{AK} guarantees that we can take the double of the pencil $(X,f)$ provided there is \textit{some} pencil with the same genus and same number of base points for which the doubling procedure is known to work. On the other hand, as observed in \cite{Smith}, for any $g\geq 2$ and $m \leq 2g-2$, there is a holomorphic genus-$g$ pencil with $m$ base points on a blow-up of the complex $\K$ surface, which can be doubled. Modeling the doubling of the convex piece after this complex model, we can thus glue the two pieces symplectically to get back $X$ with a new pencil $f'$ as above. As the gluing involves scaling the form on one of the pieces, it is clear from the construction that we can equip the two pencils with symplectic forms $\omega$ and $\omega'$ that are at least deformation equivalent. We therefore note the following as a consequence of the works of Smith and Auroux-Katzarkov:\footnote{I am grateful to Denis Auroux for verifying the arguments given here.}

\begin{lemma}\label{doubling}
Let $(X,\omega,f)$ be a genus-$g$ symplectic Lefschetz pencil with $m$ base points. If $g \geq 2$ and $m \leq 2g-2$, then we can take its double to obtain a genus-$g'$ symplectic Lefschetz pencil $(X, \omega', f')$ with $m'$ base points, where $g=2g'+m-1$ and $m'= 4m$, and $\omega$ and $\omega'$ are deformation equivalent. 
\end{lemma}

For a symplectic pencil $(X,\omega, f)$, fibers are $J$-holomorphic with respect to a suitably chosen almost complex structure $J$ compatible with $\omega$. It follows from Taubes' correspondence between Gromov and Seiberg-Witten invariants on symplectic $4$-manifolds with $b^+(X)>1$ that exceptional classes $e_j$ in $H_2(X)$ are represented by disjoint  $J$-holomorphic $(-1)$-spheres $S_j$ \cite{Ta}. The same holds when $b^+(X)=1$ and $X$ is not a rational or ruled surface by the work of Li and Liu \cite{Li1}. From the positivity of intersections for $J$-holomorphic curves, we conclude that each $S_j$ is a degree $s_j$ multisection (which we will call an \textit{$s_j$-section} in short), intersecting genus $g \geq 2$ generic fiber $F$ positively at exactly $s_j= S \cdot F \geq 1$ points. Moreover, in this case, \, $\sum s_j = (\sum S_j) \cdot F \leq 2g-2$ \, by the Seiberg-Witten adjunction inequality. Since we can always equip a Lefschetz fibration with a symplectic form with respect to which any given finite collection of disjoint sections are symplectic, we note that when $X$ is not a rational or ruled surface, there are at most $2g-2$ base points for a genus-$g$ Lefschetz pencil on $X$.

We are interested in tracking how exceptional classes, as multisections, intersect the fiber $F'$ of the new pencil produced after partial doubling. As $[F']=2[F]$ in $H_2(X)$, any exceptional sphere $S$ that is an $s$-section of $(X, \omega, f)$ gives rise to a \linebreak $2s$-section of $(X, \omega, f')$. Note that $S$ misses the base points in $(X,\omega, f)$. We introduce the following notation for the combinatorial data encoding the \emph{number} of certain exceptional classes. Consider the infinite tuple of integers 
\[(m_0, m_1, \ldots, m_r, 0, 0, \ldots) = (m_0, m_1, \ldots, m_r) \, , \]
where $m_r$ is the rightmost non-zero entry, thus allowing us to truncate the infinite tuple as we did on the right hand side. Letting $m_0$ denote the number of base points in $(X, \omega, f)$ and $m_{i+1}$, for $i \geq 0$, denote the number of $2^{i}$-sections of it representing exceptional classes, we will call the above tuple \textit{exceptional data} for $(X, \omega, f)$. Here $(X, \omega, f)$ can of course have other exceptional $s$-sections for $s \neq 2^{i}$, but for our purposes, it will suffice to keep track of the above partial information alone ---which will become evident shortly.

Now, we can partially double $(X,\omega, f)$ at $m_0 \geq k \geq 1$ points to obtain a new Lefschetz pencil $(X', \omega', f')$ (where $(X',\omega')=(X,\omega)$ if $k=m_0$). First blowing-up $X$ at $m_0-k$ base points, we get a Lefschetz pencil with the exceptional data $(k, m_0+m_1-k, \ldots)$. Then doubling at the remaining $k$ points, we arrive at the exceptional data 
\[ (4k, 0, m_0+m_1-k, \ldots, m_r) \,  \]
for the pencil $(X', \omega', f')$. Note that the length of the truncated tuple is increased by one. Importantly, if the exceptional sphere $S$ is an $s$-section of $(X, \omega, f)$ with $s \neq 2^i$, then it becomes a $2s$-section of $(X', \omega', f')$ which is not a $2^i$-section either, for any $i \geq 0$. Moreover, it is clear that the exceptional data is preserved under any orientation-preserving diffeomorphism of the pair $(X, F)$, where $F$ is the fiber of the pencil, so in particular, it is an invariant of the isomorphism class of the pencil.

A final ingredient we need is due to Welschinger: If $L$ is a Lagrangian torus in $(X, \omega)$, then any exceptional class can be represented by a symplectic $(-1)$-sphere $S$ disjoint from $L$ by \cite[Theorem~1.3]{Welschinger}. Moreover, Welschinger's proof can be applied simultaneously to a collection of exceptional spheres $S_j$ so as to symplectically isotope them away from $L$. Note that these $S_j$ are not guaranteed to be $J$-holomorphic for an almost complex structure $J$ making the fibration $f$ \linebreak $J$-holomorphic. In particular, we do not claim that they are multisections any more. The upshot here is that any (Luttinger) surgery along $L$, which would be disjoint from some regular fiber $F$, does not change the intersection number of $S$ with $F$ --identified with their inclusions-- in the new symplectic manifold. (Having checked out this aspect, we can again find homologous $J$-holomorphic representatives for $S_j$ to continue any blow-up and partial doubling process.)

We summarize what we have so far in the following lemma:

\begin{lemma}\label{partial}
Let $(X, \omega, f)$ be a symplectic Lefschetz pencil, where $X$ is not a rational or ruled surface. The exceptional data $(m_{0} , m_1, \ldots, m_r)$ is an invariant of the isomorphism class of the pencil, and moreover, it is invariant under Luttinger surgeries. For $m_0 \geq k \geq 1$, the exceptional data of a partial double of $(X, \omega, f)$ along $k$ base points is uniquely determined as $(4k, 0, m_{0}+m_1-k, \ldots, m_r)$. 
\end{lemma}

Since any sequence of partial doublings will result in a pencil on some blow-up of $(X, \omega)$, with Lemma~\ref{partial} in hand, our proof of Theorem~A now boils down to finding distinct doubling sequences resulting in pencils which
\begin{enumerate}
\item land on the same (symplectic) manifold $X'$, and \
\item have the same number of base points, \
\item are of the same genera, and
\item have distinct exceptional data.
\end{enumerate}
To strike $(1)$ and $(2)$ we will simply look at Lefschetz \textit{fibrations} with symplectic $(-1)$-sphere sections obtained by blowing-up all the base points of the pencils obtained from respective partial doubling sequences. We can then blow-down any number of these base points to produce the desired pencils in the statement of our theorem. Note that blow-up/blow-down process simply shifts weights back and forth between the first two entries of the exceptional data.  Once we generate pencils on the same manifold, we will not worry about this when comparing the exceptional data, as we will generate examples that already differ in their further entries. 

If $(X',\omega',f')$ is obtained from $(X, \omega, f)$ by a sequence of partial doublings, we see that both the smooth $4$-manifold $X'$ and the exceptional data for the pencil $f'$ are \emph{uniquely} determined by the initial exceptional data of $(X, \omega, f)$ and the ordered tuple of integers $k_1, \ldots, k_d$, for each partial doubling along $k_j$ points. Let us denote the latter sequence by $D=[k_1, \ldots, k_d]$, which is subject to the condition $4k_{j} \geq k_{j+1} \geq 1$ for all $j$. Here we use brackets both to bear in mind this extra condition, as well as to distinguish it from our notation for the exceptional data $(m_0, m_1, \ldots, m_r)$. 

Recall that each time we take a partial double of a pencil with $m$ base points along $k$ points, we pass to a pencil with $4k$ base points on a symplectic manifold which is $(m-k)$ times blow-up of the original one. By induction, we conclude that a partial doubling sequence $D=[k_1, \ldots, k_d]$ applied to a pencil with $m$ base points results in a Lefschetz fibration on a symplectic manifold which is
\[ (\cdots(((m-k_1) + 4k_1)-k_2)+ 4 k_2 + \ldots )\cdots)+ 4k_d = m + 3 \sum_{i=1}^d k_i \]
times blow-up of the initial manifold. On the other hand, if we start with a genus $g_0$ pencil, the genus of the resulting pencil after applying $D$ inductively is
\[ (\cdots(2(2g_0 +k_1-1)+k_2-1)\cdots)+ k_d -1 = 2^d g_0 + \sum_{i=1}^d 2^{d-i} (k_i-1) . \]

We therefore have:

\begin{lemma} \label{matching}
Let $f$ and $f'$ be genus $g_0$ and $g'_0$ Lefschetz pencils on $(X,\omega)$ with $m_0$ and $m'_0$ base points. Two partial doubling sequences 
\[ D=[k_1, \ldots, k_d] \, \, \text{and} \, \, D'=[k'_1, \ldots, k'_{d'}] \, \]
applied to $f$ and $f'$ (whenever it is possible, in particular when $g_0, g'_0 \geq 2$ and $m_0 \leq 2g-2$, $m'_0 \leq 2g'-2$), respectively, result in Lefschetz fibrations on the same blow-up $X \# M \CPb $ and with the same fiber genus $g$ if and only if
\vspace{-0.5cm} 
\begin{center}
\begin{align*}
M &= m_0+ 3\sum_{i=1}^d k_i = m'_0+ 3\sum_{i=1}^{d'} k'_i  \, \, \, \, \text{and} \\ 
g &= 2^d g_0 + \sum_{i=1}^d 2^{d-i} (k_i-1)= 2^{d'} g'_0 + \sum_{i=1}^{d'} 2^{d'-i} (k'_i-1) \, . 
\end{align*}
\end{center}
\end{lemma}

The proof of our theorem now reduces to generating the desired partial doubling sequences. We will show that there exists even in a simpler setting than what we have in Lemma~\ref{matching}, namely, when $f=f'$ and $d=d'$. That is, we will begin with a pencil $(X, \omega, f)$ and apply partial doubling sequences of the same length to achieve all the conditions (1)--(4) listed above. In this case, we need 
\[ \sum_{i=1}^d (k_i-k'_i) = \sum_{i=1}^d 2^{d-i} (k_i-k'_i)= 0 \, . \]
Moreover, since we depart from the same pencil and apply partial doubling sequences of the same length, without loss of generality we can take the initial exceptional data as $(m_0, 0, \ldots) = (m_0)$. In particular, it is enough to run our arguments for a minimal symplectic $4$-manifold.

First, we note that we can take $m_0$ arbitrarily large, since not only is there a symplectic pencil on any symplectic $4$-manifold, but also that one is guaranteed to have one with arbitrarily large number of base points by increasing the degree in Donaldson's construction  \cite{Donaldson}, say by further full doublings we apply to the initial pencil. Also note that neither blow-ups nor doublings would produce new separating vanishing cycles, so the second assertion in our theorem will come for free.


Let $D(n)=[k_1(n), k_2(n), k_3(n)]$ be a family of partial doubling sequences with $4k_{j}(n) \geq k_{j+1}(n) \geq 1$  and $m_0 \geq k_1$, for all $j=1,2,3$ and $n \in \N$. Regarding $D(n)$ as a $3$-dimensional integral vector, it  suffices to show that the kernel of the integral matrix 
\vspace{-0.5cm}
\begin{center}
\begin{math}
\bordermatrix{&   &    &    \cr 
              & 1 & 1  & 1  \cr
              & 2^2 & 2^1 & 1 \cr}
\end{math}
\end{center}
contains $D(n+1)-D(n)$ for at least $N$ consecutive values of $n$, provided $m_0$ is large enough. Observe that the vector $[1,-3,2]$ lies in the kernel. If we set 
\[ D(n)=[m_0-n, m_0 +3n, m_0 -2n] \, , \]
we get a legitimate partial doubling sequence, if
\[4m_0 - 4n \geq m_0 + 3n \geq 1 \, \, \, \text{and} \, \,  4m_0 + 12 n \geq m _0 -2n \geq 1 \, . \]
Recalling that $m_0 \geq 1$, we see that all needed here is 
\[m_0  \geq  \frac{7}{3}n \, \, \, \text{and} \, \, m _0 \geq 1 +2n\, , \]
which is easily seen to be satisfied by $N$ many consecutive non-negative integer values of $n$ once $m_0$ is large enough. The exceptional data corresponding to the final pencil is 
\[(4 m_0 -8n, 0, 3m_0 +14n, 3m_0-7n, n)\]
by Lemma~\ref{partial}. This completes the proof of Theorem~A.

\qed

\vspace{0.3in}
\begin{remark}
We shall note that, although we have presented our arguments in terms of geometric representatives (as multisections representing exceptional classes), the inequivalent Lefschetz pencils we have constructed are at the end distinguished by the homology classes of their fibers.
(In fact, as suggested by the referee, in the case where $X$ is an algebraic surface, one would hope to obtain similar results by showing that certain suitably chosen homology classes on a suitable blow-up are very ample and have different intersection numbers with the exceptional classes.) We therefore cannot push the same idea any further to generate infinitely many such pencils. On the other hand, one can also ask how many distinct pencils \textit{with the same fiber class} a symplectic $4$-manifold $(X, \omega)$ can be equipped with. In \cite{BH}[Theorem~1.4], using completely different methods, we present such examples which even fix the homeomorphism class of the the pair $(X,F)$, for $F$ the fiber. 
\end{remark}

\begin{remark}
Lefschetz fibrations are seen to convey different features depending on the symplectic Kodaira dimension of the underlying symplectic $4$-manifold \cite{B1, BH}. Aforementioned examples of inequivalent fibrations of Park-Yun and Smith were on symplectic $4$-manifolds of Kodaira dimension $1$. Since symplectic $4$-manifolds of negative Kodaira dimension are precisely the rational and ruled surfaces \cite{Li3}, our theorem presents inequivalent fibrations on symplectic $4$-manifolds of all non-negative Kodaira dimensions. It is plausible that, with some extra care, our construction can be carried out on rational and ruled surfaces as well. In this case, the main complication we have is that the exceptional $(-1)$-multisections might intersect each other, which in turn alters how the exceptional data changes under partial doublings involving blow-ups. It would be interesting to determine if the genus $g$ Lefschetz fibrations on ruled surfaces $\Sigma_{g/2} \x S^2 \# 4 \CPb$ and $\Sigma_{(g-1)/2} \x S^2 \# 8 \CPb$, which are known to realize the minimum number of Lefschetz singularities \cite{St2}, are unique up to isomorphisms and fibered Luttinger surgeries.
\end{remark}

\vspace{0.1in}

\section{Partial conjugations and cobordisms} \label{Partial}

We now prove Theorem~B on inequivalent Dehn twists factorizations in the mapping class group, in connection with the question on the Luttinger surgery equivalence of symplectic $4$-manifolds.

\begin{proof}[Proof of Theorem~B]
Let $W, W', W''$ be nontrivial products of positive Dehn twists along nonseparating curves, where $W= W' W''$. If the product $W' = \prod t_{a_i}$, as a mapping class, stabilizes a loop $\alpha$ on $\Sigma_g^m$, then $W'$ commutes with the Dehn twist $t_\alpha \in \Gamma_g^m$.  We can then produce a new positive factorization $W'_{\alpha} W''$, which is derived from $W$ by an \textit{untwisted partial conjugation}, for $W'_{\alpha} = \prod  t_{t_\alpha(a_i)} = \prod t^{-1}_{\alpha} t_{a_i} t_{\alpha}$. 

Since any boundary twist commutes with nonseparating Dehn twists on the surface, we can assume without loss of generality that $\alpha$ is a --possibly separating-- curve which is not boundary parallel. Taking a parallel transport of $\alpha$ over a curve $\gamma$ enclosing the Lefschetz critical values corresponding to all vanishing cycles $a_1, \ldots a_k$ in the product $\prod t_{a_i}$, we produce a Lagrangian torus fibered over $\gamma$. The untwisted partial conjugation amounts to a Luttinger surgery along the torus in the direction of $\alpha$ \cite{Auroux2}. 

Now, let $f_i$, $i=1,\ldots, N$ be any $N$ distinct genus $g$ Lefschetz pencils on a symplectic $4$-manifold $X$ provided by our Theorem~A, with the additional feature that the vanishing cycles are nonseparating. Observe that, if needed, we can take full doubles of all $f_i$ simultaneously to obtain at least $m$ base points, for any given $m$. Hence, for each $f_i$, we get a \textit{positive factorization} $W_i$ of the boundary parallel multitwist \,$t_{\delta_1} \cdot \ldots \cdot t_{\delta_m}$\, into positive Dehn twists along nonseparating curves in $\Gamma_g^m$. It immediately follows that $W_i$, $W_j$ are not equivalent via Hurwitz moves and untwisted partial conjugations by Dehn twists for any $i\neq j$. 
\end{proof}

\vspace{0.6cm}
\begin{remark}
One can similarly define a surgery along a Klein bottle $L$ as an equivariant surgery along its double cover, as well as its symplectic analogue as a Luttinger surgery; see e.g. \cite{Nemirovski}.  If the product $W'$ in the above proof, as a mapping class, maps the loop $\alpha$ to $-\alpha$, then $W'$ still commutes with the Dehn twist $t_\alpha \in \Gamma_g^m$, and one can produce a new positive factorization $W'_{\alpha} W''$, which is now derived from $W$ by a \textit{twisted partial conjugation}. The twisted partial conjugation in this case amounts to a fibered Luttinger surgery along a Lagrangian Klein bottle $L$ \cite{Shev}.  It is plausible that Theorem~B (and Theorem~A) can be extended to include this case as well. However, our earlier arguments do not go through to conclude this, since there is no analogue of Welschinger's result we used in the proof of Theorem~A to guarantee that exceptional spheres can be isotoped away from a given Lagrangian Klein bottle. (In fact this is not true in general: if we take the obvious Lagrangian Klein bottle $L$ in $X= \CP \# \CPb$, then there exists a pencil on $X$ for which $L$ is fibered. If $L$ were disjoint from the exceptional sphere $S$ in $X$, we could blow-down $S$ to get a Lagrangian embedding of $L$ in $\CP$, which contradicts the main theorem of  \cite{Nemirovski, Shev}.)
\end{remark}

\begin{remark}
Let $[W',\phi]=1$ in $\Gamma_g^m$. Although any mapping class $\phi$ in $\Gamma_g^m$ can be written as a product of Dehn twists, one cannot necessarily choose these Dehn twists in a way that \textit{each} one of them commutes with $W'$: For instance, if $\phi(a)=-a$, we certainly have $[t_a, \phi]=1$. On the other hand, for any Dehn twist $t_b$ satisfying $[t_a, t_b]=1$, we have $b$ isotopic to a curve disjoint from $a$. However, a product of such twists could not reverse the orientation on $a$. Hence, the equivalence of two factorizations via partial conjugations is more general than the equivalence of them via partial conjugations by Dehn twists. 
\end{remark}


\begin{remark}
Any two integral symplectic $4$-manifolds with the same characteristic numbers $(c_1^2, c_2, [\omega]^2$, $c_1 \cdot [\omega])$ have the same Euler characteristic and signature, determined by $c_1^2$ and $c_2$. Therefore the results of \cite{BS} show that if one only considers the underlying smooth structures, the two $4$-manifolds would be equivalent via smooth surgeries along tori. However, when we in addition take compatible symplectic Lefschetz pencils on them, Theorem~B dictates that there is no sequence of Luttinger surgeries taking one to the other. Curiously, Auroux's other question on the surgery equivalence of such integral symplectic $4$-manifolds \cite{Auroux1} lies in the middle ground these results fall short of covering. 
\end{remark}


\vspace{0.3in}
\begin{thebibliography}{99}

\bibitem{ADK} D. Auroux, S. K. Donaldson, L. Katzarkov, {\it Luttinger surgery along Lagrangian tori and non-isotopy for singular symplectic plane curves,} Math. Ann. 326 (2003), 185--203.

\bibitem{Auroux1} D. Auroux, {\it Symplectic $4$-manifolds, singular plane curves, and isotopy problems,} Floer homology, gauge theory, and low-dimensional topology, 263--276, Clay Math. Proc., 5, Amer. Math. Soc., Providence, RI, 2006.

\bibitem{Auroux2} D. Auroux, {\it Mapping class group factorizations and symplectic $4$-manifolds: some open problems,} Problems on mapping class groups and related topics, 123--132, Proc. Sympos. Pure Math. 74, Amer. Math. Soc., Providence, RI, 2006.

\bibitem{AK} D. Auroux and L. Katzarkov, {\it A degree doubling formula for braid monodromies and Lefschetz pencils,} Pure Appl. Math. Q. 4 (2008), no. 2, part 1, 237--318.

\bibitem{B1} R.\,I. Baykur, {\it Minimality and fiber sum decompositions of Lefschetz fibrations,} to appear in Proc., Amer. Math. Soc.; http://dx.doi.org/10.1090/proc/12835.

\bibitem{BH} R.\,I. Baykur and K. Hayano, {\it Multisections of Lefschetz fibrations and topology of symplectic $4$-manifolds, } preprint; http://arxiv.org/abs/1309.2667. 

\bibitem{BS} R.\,I. Baykur and N. Sunukjian, {\it Round handles, logarithmic transforms and smooth \linebreak $4$-manifolds,} J. Topol. 6 (2013), no. 1, 49--63. 


\bibitem{Donaldson} S. K. Donaldson, {\it Lefschetz pencils on symplectic manifolds}, J. Differential Geom. 53 (1999), no.2, 205--236


\bibitem{GS} R. E. Gompf and A.I.Stipsicz, {\it 4-Manifolds and Kirby Calculus}, Graduate Studies in Mathematics 20, American Mathematical Society, 1999.

\bibitem{HL} C.-I. Ho and T.-J. Li, {\it Luttinger surgery and Kodaira dimension},  Asian J. Math. 16 (2012), no. 2, 299–-318. 


\bibitem{Kas} A. Kas, {\it On the handlebody decomposition associated to a Lefschetz fibration}, Pacific J. Math. 89 (1980), 89--104.


\bibitem{Li1} T.-,J. Li and A. Liu, {\it Symplectic structures on ruled surfaces and a generalized adjunction formula,} Math. Res. Lett. 2 (1995), 453--471.

\bibitem{Li3} T.-J. Li, {\it The Kodaira dimension of symplectic $4$-manifolds,} Floer homology, gauge theory, and low-dimensional topology, 249--261, Clay Math. Proc., 5, Amer. Math. Soc., Providence, RI, 2006.

\bibitem{Luttinger} K. M. Luttinger, {\it Lagrangian tori in $\R^4$}, J. Differential Geom. 42,  (1995), no. 2, 220--228.  

\bibitem{Matsumoto} Y. Matsumoto, {\it Lefschetz fibrations of genus two - a topological approach} -, Proceedings of the 37th Taniguchi Symposium on Topology and Teichm\"{u}ller Spaces, (S. Kojima, et. al., eds.), World Scientific, 1996, 123--148

\bibitem{Nemirovski} S. Nemirovski, {\it The homology class of a Lagrangian Klein bottle,} Izv. Ross. Akad. Nauk Ser. Mat. 73 (2009), no. 4, 37--48.



\bibitem{ParkYun} J. Park and K.-H. Yun, {\it Nonisomorphic Lefschetz fibrations on knot surgery $4$-manifolds,} Math. Ann. 345 (2009), no. 3, 581--597.

\bibitem{Shev} V. Shevchishin, {\it Lagrangian embeddings of the Klein bottle and the combinatorial properties of mapping class groups,} Izv. Ross. Akad. Nauk Ser. Mat. 73 (2009), no. 4, 153--224. 

\bibitem{Smith} I. Smith, {\it Lefschetz pencils and divisors in moduli space,} Geom. Topol. 5 (2001), 579--608.

\bibitem{Stern} R. Stern, {\it Will we ever classify simply-connected smooth $4$-manifolds?}, Floer homology, gauge theory, and low-dimensional topology, Clay Mathematics Proceedings: 5, American Mathematical Society, Providence, RI, 2006, 225--239.

\bibitem{St2} A. Stipsicz, {\it Singular fibres in Lefschetz fibrations on manifolds with $b_2^+=1$}, Topology Appl. 117 (2002), no. 1, 9--21.

\bibitem{Ta} C. H. Taubes, {\it The Seiberg--Witten and Gromov invariants,} Math. Res. Lett. 2 (1995), no. 2, 221--238.

\bibitem{Welschinger} J. Welschinger, {\it Effective classes and Lagrangian tori in symplectic four-manifolds,} J. Symplectic Geom. 5 (2007), no. 1, 9--18. 





\end{thebibliography}
\end{document}